\documentclass[11pt]{article}

\usepackage{amsmath,amssymb,amsthm,epsfig}

\newcommand{\vp}{\varepsilon}
\newcommand{\n}{\noindent}
        
\newcommand{\ovl}{\overline}
\newcommand{\bb}[1]{\mathbb{#1}}
\newcommand{\cl}[1]{\mathcal{#1}}
\newcommand{\intl}{\int\limits}

\theoremstyle{plain}
\newtheorem*{thmb}{Theorem}

\theoremstyle{remark}
\newtheorem{rem}{Remark}

\theoremstyle{corollary}
\newtheorem{cor}{Corollary}

\theoremstyle{proposition}

\begin{document}

\title{Geometric conditions which imply compactness of the 
$\bar\partial$-Neumann operator\thanks{Research supported in
part by NSF grant DMS 0100517}}

\author{Emil J.\ Straube \\ Department of Mathematics, Texas
A\&M University, College Station,TX
77843\\straube@math.tamu.edu}

\date{}

\maketitle

\begin{abstract}
For smooth bounded pseudoconvex domains in $\bb C^2$, we
provide geometric conditions on (the points of infinite type
in) the boundary which imply compactness of the
$\bar\partial$-Neumann operator. It is noteworthy that the
proof of compactness does {\it not} proceed via verifying
the known potential theoretic sufficient conditions.
\end{abstract}

\bigskip
\bigskip
Let $\Omega$ be a smooth bounded pseudoconvex domain in $\bb C^n$. The 
$\bar\partial$-Neumann operator $N$ on (0,1)-forms is the inverse of the complex 
Laplacian $\bar\partial\bar\partial^* +
\bar\partial^*\bar\partial$. $N$ and its regularity theory
play a crucial role both in partial differential equations
and in several complex variables ([\ref{FK}], [\ref{BS}],
[\ref{CS}]). In particular, the question when $N$ is a
compact operator on $\cl L^2_{(0,1)}(\Omega)$, the space of
(0,1)-forms with square integrable coefficients, is of
interest in several contexts. Well known examples are global
regularity and the Fredholm theory of Toeplitz operators. We
refer the reader to [\ref{FS1}] for more information and for
references on these and other questions related to
compactness of the $\overline{\partial}$-Neumann operator.
Among more recent references, we mention [\ref{HL}], where
compactness is related to the existence of Henkin--Ramirez
type integral formulas with well behaved kernels. And, there
is a useful connection between compactness in the
$\overline{\partial}$-Neumann problem and the asymptotic
behavior, in a semi--classical limit, of the ground state
energy of certain magnetic Schr\"odinger operators
([\ref{FS3}], [\ref{CF}]).

The most general known sufficient condition for compactness
is potential theoretic in nature. Roughly speaking, there
should exist, near the boundary points of infinite type,
plurisubharmonic functions with arbitrarily large complex
Hessians whose gradients are uniformly bounded in the metric
induced by the Hessians of the functions. This condition was
introduced under the name property $(\widetilde P)$ and
shown to imply compactness by McNeal in [\ref{McN}].
It generalizes previous work by Catlin that used his now
classical condition property $(P)$ ([\ref{C1}], see
[\ref{S1}] for a systematic study of this property).
Property $(\widetilde P)$ implies compactness of the
$\overline{\partial}$-Neumann operator on an arbitrary
bounded pseudoconvex domain, even when no boundary
regularity whatsoever is assumed ([\ref{St}], Corollary 3;
[\ref{McN}], Corollary 4.2).

There are natural geometric conditions,
bearing on how the set of infinite type boundary points
sits inside the boundary, which are known to imply property
$(P)$ or $(\widetilde P)$ and hence compactness of the
$\overline{\partial}$-Neumann operator ([\ref{C1}],
[\ref{S1}], [\ref{B}], [\ref{FS1}], [\ref{FS2}],
[\ref{McN}]). In this paper, we present new geometric
conditions which imply compactness in the case of
smooth bounded pseudoconvex
domains in $\mathbb{C}^{2}$. In light of the
discussion in the next paragraph, the fact that compactness is {\it not}
demonstrated via establishing property $(P)$ or $(\widetilde
P)$ (in contrast
to ([\ref{C1}], [\ref{S1}], [\ref{B}], [\ref{FS1}],
[\ref{FS2}], [\ref{McN}]) is worth pointing out.

On domains which are locally convexifiable, the analysis,
the potential theory, and the geometry associated with the
$\overline{\partial}$-Neumann problem mesh perfectly. That
is, the following three statements are equivalent: (i) the
$\overline{\partial}$-Neumann operator is compact, (ii) the
boundary of the domain satisfies property $(\widetilde P)$,
(iii) the boundary of the domain does not contain (germs of)
analytic discs ([\ref{FS2}], [\ref{FS1}], Theorem 5.1). This
is far from true in general. Property$(\widetilde P)$ always
excludes discs from the boundary (as is  seen by pulling
back the good plurisubharmonic functions to the unit disc in
the complex plane); so does compactness of $N$ for domains
in $\mathbb{C}^{2}$ with at least Lipschitz boundary (see
[\ref{FS1}], Proposition 4.1 for a proof of this folklore
result). However, Sibony ([\ref{S1}]) observed that the
absence of discs from the boundary need not imply property
$(\widetilde P)$, and Matheos ([\ref{M}], see also [\ref{FS1}])
showed that in fact the absence of
discs need not (even) imply compactness of the
$\overline{\partial}$-Neumann operator. (Actually, Sibony's
observation concerned property $(P)$, but the
domains in [\ref{S1}] are Hartogs domains in
$\mathbb{C}^{2}$ where $(\widetilde P)$ and $(P)$ are
equivalent ([\ref{FS3}], Appendix A).) This left open the
exact relationship between property $(\widetilde P)$ and
compactness of the $\overline{\partial}$-Neumann operator.
Christ and Fu ([\ref{CF}]) recently established the
equivalence of these two properties on smooth bounded
pseudoconvex Hartogs domains in $\mathbb{C}^{2}$. On general
domains, however, the situation is not understood. In view
of this, it is desirable to have a technique like ours
for establishing compactness of the
$\overline{\partial}$-Neumann operator that does
not rely on property $(\widetilde P)$ (in this context,
compare also [\ref{HL}]).

\bigskip

If $Z$ is a (real) vector field defined on some open subset
of $b\Omega$ (or of $\bb C^2$), we denote by $\cl F^t_Z$ the
flow generated by $Z$. By finite type
of a boundary point, we mean finite type in the sense of D'Angelo; however, for 
domains in $\bb C^2$, the various notions of finite type all agree, so that no 
distinction is necessary ([\ref{DA}]). Recall that the set
of points of finite type in the boundary of a smooth bounded
pseudoconvex domain is open and, consequently, the set of
points of infinite type is compact. $B(P,r)$ denotes the
open ball of radius $r$ centered at $P$.

\begin{thmb}
Let $\Omega$ be a $C^\infty$-smooth bounded pseudoconvex domain in $\bb C^2$.
Denote by $K$ the  set of boundary points of infinite type. Assume that
there exist constants $C_1,C_2>0$, $C_3$ with $1\le
C_3<3/2$, and a sequence $\{\vp_j > 0 \}^\infty_{j=1}$ with
$\lim\limits_{j\to\infty} \vp_j = 0$ so that the following
holds. For every $j\in\bb N$ and $P\in K$ there is a (real)
complex tangential vector field $Z_{P,j}$ of unit length
defined in some neighborhood of $P$ in $b\Omega$ with
$\max|\hbox{div } Z_{P,j}|\le C_1$ such that $\cl F
^{\vp_j}_{Z_{P,j}}(B(P,C_2(\vp_j)^{C_3})\cap K) \subseteq
b\Omega\backslash K$. Then the $\overline{\partial}$-Neumann
operator on $\Omega$ is compact.
 \end{thmb}

It is of course assumed that the flow $\cl F^t_{Z_{P,j}}$
exists for all initial points in $B(P,C_2(\vp_j)^{C_3})\cap
K$ and $t\le \vp_j$.

\bigskip

Two immediate questions arise. First, the obvious
examples that satisfy the assumptions in the theorem also
satisfy property$(\widetilde P)$, and whether or not the
theorem can actually furnish examples (if they exist) of
domains with compact $\overline{\partial}$-Neumann operator,
but without property $(\widetilde P)$, does not appear to be
a simple matter to decide. But we obtain, in any event, a
simple geometric proof of compactness under the
conditions in the theorem. Moreover, these conditions
turn out to be natural and, modulo the form of the
lower bound $C_2(\vp_j)^{C_3}$, minimal; we discuss
this in Remark 3 below. Second, our result is stated for
domains in $\bb C^2$; the proof uses maximal estimates and
so does not work in higher dimensions. However, it is easy
to see that the theorem, when transcribed to $\bb C^n$
verbatim, also fails. Control over $K$ will have to be
imposed in all complex tangential directions rather than
just in one. We hope to return to these questions elsewhere.

\medskip

The general thrust in the theorem is that at points of $K$
there should exist a (real) complex tangential direction
transversal to $K$ in which $b\Omega \backslash K$ (the
good set) is thick enough. This occurs in extremis
when $K$ is (locally) contained in a totally real
submanifold of the boundary (property$(P)$, and hence
compactness, are well known in this case ([\ref{C1}])).
The following corollary to the theorem takes this
situation and disposes of the requirement that $K$ be
(contained in) a smooth submanifold of the boundary.

\begin{cor}\label{cor1}
Let $\Omega$ be a $C^\infty$-smooth bounded pseudoconvex domain in $\bb C^2$.
Denote by $K$ the  set of boundary points of infinite type.
Assume that for all $P\in K$ there exists a (real)
complex tangential vector field $Z_p$ defined near $P$, a
neighborhood $U_p$ of $P$, and $\vp_p>0$ such that $\cl F
^t_{Z_p}(U_p\cap K) \subseteq b\Omega\backslash K$ for $t\le
\vp_p$. Then the $\overline{\partial}$-Neumann operator on
$\Omega$ is compact.
\end{cor}

To reduce Corollary 1 to the theorem, it suffices to cover
$K$ by finitely many neighborhoods $U_P$.

A second
geometrically simple special case occurs when $b\Omega
\backslash K$ satisfies what might be called a weak complex
tangential cone condition. That is, there should exist a
finite (possibly small) open real cone $C$ in $\bb
C^2\approx \bb R^4$ having the following property. For each
$P\in K$ there exist a complex tangential direction so that
when $C$ is moved by a rigid motion to have vertex at $P$
and axis in that complex tangential direction, the (open)
cone obtained intersects $b\Omega$ in a set contained in
$b\Omega\backslash K$. Then the assumptions in the theorem
are satisfied with $C_3=1$.

\begin{cor}\label{cor2}
Let $\Omega$ be a $C^\infty$-smooth bounded pseudoconvex domain in $\bb C^2$.
Denote by $K$ the  set of boundary points of infinite type.
Assume that $b\Omega\backslash K$ satisfies a weak complex
tangential cone condition. Then the
$\overline{\partial}$-Neumann operator on $\Omega$ is
compact.
\end{cor}

Because $C_3$ is only required to be in the range $1 \leq
C_3 < 3/2$, the cone in Corollary 2 may be allowed to
degenerate into a mild cusp; we leave the details to the
reader.

\medskip

\begin{rem}\label{rem1}
The assumption in the theorem says that for all $j\in \bb N$, all $P\in K$,
there should exist a complex tangential vector field and a ball centered at $P$ 
so that the intersection of $K$ with this ball, when translated along the 
integral curves of the vector field by $\vp_j$, ends up in $b\Omega\backslash 
K$, 
i.e.\ in the set of points of finite type. The crucial additional condition is 
that the radius of this ball should be at least of order $(\vp_j)^{C_3}$ for 
some $C_3<3/2$. A condition of this kind is needed; it is not enough to just 
have some ball centered at $P$. Consider a Hartogs domain $\Omega := \{(z,w)\in 
\bb C^2/|z|<1, |w|<e^{-\varphi(z)}\}$,  where $\varphi$ is smooth and 
plurisubharmonic on the unit disc, and such that $\Omega$ is a smooth domain in 
$\bb C^2$, which is strictly pseudoconvex at boundary points where $w=0$. Let 
$K_0 := \{z\in D\mid \Delta\varphi = 0\}$ be relatively compact in $D$. Assume 
that $\Delta\varphi$ vanishes to infinite order on $K_0$. It is well known that 
the set $K$ of boundary points of infinite type consists precisely of those 
$(z,w)\in b\Omega$ with $z\in K _0$. Finally, assume that $K_0$ has empty 
interior (this is equivalent to $b\Omega$ not containing analytic discs). If $P 
= (z,w)\in K \subseteq b\Omega$, then $z\in K_0$. Because $K_0$ has empty 
interior, there are arbitrarily short translates of $z$ in $D\backslash K_0$. 
Because $D\backslash K_0$ is open, there exists for every $\vp>0$ a vector $X$ 
and a ball $B(z,r)$ such that $B(z+X,r) \subseteq D \backslash K_0$. The 
(constant) vector field $X$ is easily lifted to a complex tangential vector 
field on $b\Omega$ (near $P$). In this way, one sees that $\Omega$ satisfies the 
assumptions of the theorem, except for the lower bound (in terms of $\vp_j$) of 
the balls centered at points of $K$. This is enough to make the theorem fail:\ 
there are examples of Hartogs domains as above, where the $\bar\partial$-Neumann 
operator is not compact ([\ref{M}], [\ref{FS1}], Theorem 4.2 and Remark 5).
\end{rem}

\begin{rem}\label{rem2}
The key to getting the vector fields from the theorem in the
examples in Remark 1 is that $K_0$ has empty interior, or
equivalently,  that $b\Omega$ contains no analytic discs. It
turns out that suitable vector fields may be obtained on any
domain in $\bb C^2$ whose boundary contains no analytic
disc. In fact, the boundary of a smooth bounded pseudoconvex
domain in $\bb C^2$ contains no analytic discs if and only
if it satisfies all the assumptions of the theorem, except
possibly the lower bounds (in terms of $\vp_j$) on the
radii of the balls centered at points $P\in K$. It is clear
that the assumptions in the theorem, without these lower
bounds, exclude discs from the boundary. The proof of the
other direction follows easily from Proposition 3.1.12 in
([\ref{C2}]), as follows. Let $P\in b\Omega$, $Z$ a real
vector field defined on $b\Omega$ in a neighborhood of $P$
that is complex tangential, and let $J$ denote the
involution associated with the complex structure of $\bb
C^n$. Then $Z$ and $J(Z)$ span (over $\bb R$) the complex
tangent space to $b\Omega$ near $P$. For $0\le\theta <
2\pi$, denote by $Z^\theta$ the field $Z^\theta
:=(\cos\theta)Z  + (\sin(\theta) J(Z)$. Fix $\vp>0$. Catlin
shows in [\ref{C2}], Proposition 3.1.12 that if all points
of $M_\vp := \{\cl F^t_{Z^\theta}(P)\mid 0\le t<\vp, 0\le
\theta < 2\pi\}$ are weakly pseudoconvex points of
$b\Omega$, then $M_\vp$ is a (necessarily one-complex
dimensional) complex manifold. Consequently, when there is
no analytic disc in $b\Omega$, there exist complex
tangential fields $Z_{P,\vp}$ for all $P\in K$ and
sufficiently small $\vp>0$ so that $\cl F
^\vp_{Z_{P,\vp}}(z)\notin K$ for $z$ close enough to $P$.
Moreover, the $Z_{P,\vp}$ are of the form $Z_{P,\vp} =
(\cos\theta_{P,\vp}) Z + (\sin\theta_{P,\vp}) J(Z)$, so that
div $Z_{P,\vp}$ is bounded uniformly in $P$ and $\vp$.
\end{rem}

\begin{rem}\label{rem3}
Since on a smooth bounded pseudoconvex domain in
$\mathbb{C}^{2}$, compactness of the
$\overline{\partial}$-Neumann operator implies that there
are no discs in the boundary (see the discussion preceeding
the statement of the theorem), we can deduce from Remark 2 a
(very) partial converse to the theorem. If $N$ is compact,
then the domain satisfies the assumptions in the theorem
except possibly the lower bound $C_2(\vp_j)^{C_3}$ on the
radii of the balls whose intersections
with $K$, when translated along the integral curves of the
fields $Z_{P,j}$, end up in $b\Omega \backslash K$.
Moreover, the examples in Remark 1 show that there must
be some lower bound on the size of these balls in terms of
$\vp_j$, otherwise the theorem need not hold. Thus the
assumptions in the theorem are quite natural, and, modulo
the exact form of this lower bound, minimal, as asserted in
the discussion following the statement of the theorem.
\end{rem}

\bigskip

We now prove the theorem. There are various equivalent statements to compactness 
of the $\bar\partial$-Neumann operator, see e.g.\ [\ref{FS1}], Lemma 1.1. We 
will 
show a so-called compactness estimate, that is, we will show that for all 
$\vp>0$ there exists a constant $C_\vp$ such that
\begin{equation}\label{eq1}
\|u\|^2_0 \le \vp(\|\bar\partial u\|^2_0 + \|\bar\partial^*u\|^2_0) + C_\vp 
\|u\|^2_{-1},
\qquad u\in \hbox{dom}(\bar\partial^*)\cap C^\infty_{(0,1)}(\ovl\Omega).
\end{equation}
Here, $C^\infty_{(0,1)}(\ovl\Omega)$ denotes the space of (0,1)-forms with 
coefficients in $C^\infty(\ovl\Omega)$, the norm $\|~~\|_0$ and $\|~~\|_{-1}$ 
denote the Sobolev-0 (i.e.\ $\cl L^2(\Omega)$) and Sobolev-1 norms, 
respectively. For forms, these norms are computed
componentwise. Note that if \eqref{eq1} holds for
$u \in \hbox{dom}(\bar\partial^*)\cap
C^\infty_{(0,1)}(\ovl\Omega)$, it holds for $u \in
\hbox{dom}(\bar\partial^*)\cap\hbox{dom}(\bar\partial)$,
because the former space is dense in the latter with respect
to the graph norm (see for example [\ref{CS}], Lemma 4.3.2).

The idea of the proof is very simple. To estimate the $\cl L^2$-norm of $u$
near a point $P$ of $K$, we express $u$ there in terms of $u$ in a patch which 
meets the boundary in a relatively compact subset of the set of finite type 
points plus the integral of the derivative of $u$ in the direction $Z_{P,j}$.
The first contribution is easily handled by subelliptic estimates, while the
second is estimated by the length of the curve (which is $\vp_j$) times the $\cl 
L^2$-norm of $Z_{P,j}u$. But in $\bb C^2$, this $\cl L^2$-norm is estimated by 
the $\cl L^2$-norm of $\bar\partial u$ and $\bar\partial^*u$, because $Z_{P,j}$ 
is complex tangential. Controlling the terms coming from the
integral of $Z_{P,j}u$ raises  overlap and divergence
issues; these are taken care of by the uniformity built into
the assumption in the theorem.

The details are as follows. First note that we can extend
the fields $Z_{P,j}$ form $b\Omega$ to the inside of
$\Omega$ by a fixed distance by letting them be constant
along the real normal. We still denote these extended fields
by $Z_{P,j}$; they are complex tangential to the level sets
of the boundary distance.

Fix $\vp>0$ and $j$ so that $\vp_j<\vp$. A standard
covering theorem (see e.g.\ [\ref{Z}], Theorem 1.3.1)
applied to the family of closed balls
$\{\ovl{B(P,\frac{C_2}{10}(\vp_j)^{C_3})}\mid P\in K\}$
gives a subfamily of pairwise disjoint balls so that
the corresponding closed balls of radius
$\frac{C_2}{2}(\vp_j)^{C_3}$, hence the open balls of radius
$C_2(\vp_j)^{C_3}$, still cover $K$. Because $K$ is compact,
we obtain a finite family of open balls
$\{B(P_k, C_2(\vp_j)^{C_3})\mid 1\le k\le N, P_k\in K\}$
that covers $K$, and such that the corresponding closed
balls of radius $\frac{C_2}{10}(\vp_j)^{C_3}$ are pairwise
disjoint. To simplify notation, we will use $Z_k$ to denote
the fields $Z_{P_k,j}$, $1\le k\le N$. By decreasing $C_2$
in the theorem, we may assume that $\cl
F^{\vp_j}_{Z_k}(B(P_k, C_2(\vp_j)^{C_3})\cap K)$ is not only
contained in $b\Omega\backslash K$, but is relatively
compact there. Consequently there exist open subsets $U_k$,
$1\le k\le N$, of $\ovl\Omega$, with
\begin{equation}\label{eq2} K\cap B(P_k,C_2(\vp_j)^{C_3})
\subseteq U_k\subseteq B(P_k,C_2(\vp_j)^{C_3})
\end{equation} and \begin{equation}\label{eq3} \cl
F^{\vp_j}_{Z_k}(U_k\cap\Omega) \cap K = \emptyset.
\end{equation} It is implicit in the statement of
\eqref{eq3} that $U_k$ is chosen in a sufficiently small
neighborhood of $b\Omega$ so that the flow generated by the
extended vector fields $Z_k$ exists up to time $\vp_j$ for
all initial points in $U_k$.

Now let $u\in C^\infty_{(0,1)}(\ovl\Omega)\cap \text{dom}(\bar\partial^*)$. Then 
\begin{equation}\label{eq4}
\|u\|^2_0 = \int\limits_{\left(\bigcup\limits^N_{k=1} U_k\right)\cap\Omega} 
|u|^2 + \int\limits_{\Omega\big\backslash\left(\bigcup\limits^N_{k=1} 
U_K\right)} 
|u|^2.
\end{equation}
Because $\Omega\Big\backslash\left(\bigcup\limits^N_{k=1} U_k\right)$ meets 
$b\Omega$ in a compact subset of $b\Omega$  which does not intersect $K$, we can 
apply subelliptic estimates (see [\ref{C3}]) to estimate the second term in the 
right-hand side of \eqref{eq4}:\ there exist $s>0$ and $C>0$ such that the 
restriction of $u$ to a neighborhood $U$ of $\Omega \Big\backslash \left( 
\bigcup\limits^N_{k=1} U_k\right)$ (in $\Omega$) belongs to $W^s_{(0,1)}(U)$ and
\begin{equation}\label{eq5}
\|u\|^2_{W^s_{(0,1)}(U)} \le C(\|\bar\partial u\|^2_0 + \|\bar\partial^* 
u\|^2_0).
\end{equation}

\noindent The usual interpolation inequality for Sobolev 
norms gives
\begin{equation}\label{eq6}
\begin{split}
\int\limits_{\Omega\big\backslash\left(\bigcup\limits^n_{k=1} U_k\right)} |u|^2 
&\le \|u\|_{\cl L^2_{(0,1)}(U)}\\
&\le \frac\vp{C} \|u\|^2_{W^s_{(0,1)}(U)} + C_\vp\|u\|^2_{W^{-1}_{(0,1)}(U)}\\
& \\
&\le \vp(\|\bar\partial u\|^2_0 + \|\bar\partial^*u\|^2_0) + C_\vp\|u\|^2_{-1}.
\end{split}
\end{equation}

We now estimate the first term on the right-hand side of \eqref{eq4}. Fix $k$,
$1\le k\le N$.
\begin{equation}\label{eq7}
\begin{split}
\intl_{U_k\cap \Omega} |u|^2 &= \intl_{U_k\cap\Omega} \left|u(\cl F 
^{\vp_j}_{Z_k}(x)) - \intl^{\vp_j}_0 Z_ku(\cl F^t_{Z_k}(x))dt\right|^2\ 
dV(x)\\
& \\
&\le 2 \intl_{U_k\cap\Omega} |u(\cl F^{\vp_j}_{Z_k}(x))|^2 dV(x)\\
&\quad + 2 \intl_{U_k\cap\Omega} \left|\intl^{\vp_j}_0 Z_ku(\cl F^t_{Z_k}(x)) 
dt\right|^2 dV(x).
\end{split}
\end{equation}
The first term on the right-hand side of \eqref{eq7} can be estimated as
follows:
\begin{equation}\label{eq8}
\begin{split}
\intl_{U_k\cap\Omega} |u(\cl F^{\vp_j}_{Z_k}(x))|^2 dV(x) &= \intl_{\cl F
^{\vp_j}_{Z_k}(U_k\cap\Omega)} |u(y)|^2 \det(\partial
x / \partial y) dV(y)\\ &\le C_k \intl_{\cl
F^{\vp_j}_{Z_k}(U_k\cap\Omega)} |u(y)|^2 dV(y). \end{split}
\end{equation}
Here we use $\det (\partial x / \partial y)$ as
shorthand for the (positive) Jacobian of the
diffeomorphism $\cl F ^{-\vp_j}_{Z_k}: \cl
F^{\vp_j}_{Z_k}(U_k\cap\Omega) \rightarrow U_k\cap\Omega$.
By \eqref{eq3}, we can use subelliptic estimates
once more to estimate the last term in \eqref{eq8}; an
argument analogous to the one that led to \eqref{eq6} gives
\begin{equation}\label{eq9} \intl_{\cl
F^{\vp_j}_{Z_k}(U_k\cap\Omega)} |u(y)|^2 dV(y) \le
\frac\vp{C_kN} (\|\bar\partial u\|^2_0 +
\|\bar\partial^*u\|^2_0) + C_\vp\|u\|^2_{-1}. \end{equation}

We now estimate the second term on the right-hand side of \eqref{eq7}. The 
Cauchy--Schwarz inequality and Fubini's theorem give

\begin{equation}\label{eq10}
\begin{split}
\intl_{U_k\cap\Omega} \left|\intl^{\vp_j}_0 Z_ku(\cl F^t_{Z_k}(x))dt\right|^2 
dV(x) &\le \vp_j \intl_{U_k\cap\Omega} \intl^{\vp_j}_0 |Z_ku(\cl 
F^t_{Z_k}(x))|^2 dt\ dV(x)\\
&= \vp_j \intl^{\vp_j}_0 \intl_{U_k\cap\Omega} |Z_ku(\cl F^t_{Z_k}(x))|^2 
dV(x) \ dt\\
&= \vp_j \intl^{\vp_j}_0 \intl_{\cl F^t_{Z_k}(U_k\cap\Omega)} 
|Z_ku(y)|^2 
\det(\partial x / \partial y) dV(y)\ dt\\
&\le 2 \vp_j \intl^{\vp_j}_0 \intl_{\cl F^t_{Z_k}(U_k\cap \Omega)} 
|Z_ku(y)|^2 
dV(y)\ dt.
\end{split}
\end{equation}
In the last inequality in \eqref{eq10} we have used the
uniform bound on the divergence of the fields $Z_k$, that
is, on the rate of change of the volume element under the
flows generated by the $Z_k$'s. This bound implies that
$\det(\partial x / \partial y) \le \exp(tC_1)
\le \exp(\vp_j C_1) \le \exp(\vp C_1) \le 2$ for $\vp$ small
enough. (It suffices to establish \eqref{eq1} for small
enough $\vp>0$.)

Putting together estimates \eqref{eq7}--\eqref{eq10} and adding over $k$, we 
estimate the first term on the right-hand side
of\eqref{eq4}:

\begin{equation}\label{eq11}
\begin{split}
\intl_{\left(\bigcup\limits^N_{k=1}U_k\right)\cap \Omega} |u|^2 &\le
\sum^N_{k=1} \intl_{U_k\cap\Omega} |u|^2\\
&\le \sum^N_{k=1} \bigg[\frac{2\vp}N (\|\bar\partial u\|^2_0 + \|\bar\partial^* 
u\|^2_0) + C_\vp\|u\|^2_{-1}\\
&\quad + 4\vp_j \intl^{\vp_j}_0 \intl_{\cl F^t
_{Z_k}(U_k\cap \Omega)} |Z_ku(y)|^2 dV(y)\ dt\bigg]\;.\\
\end{split}
\end{equation}
That is,
\begin{equation}\label{eqinsert}
\begin{split}
\intl_{\left(\bigcup\limits^N_{k=1}U_k\right)\cap \Omega}
|u|^2 &\le 2\vp(\|\bar\partial u\|^2_0 +
\|\bar\partial^*u\|^2_0) + C_\vp \|u\|^2_{-1}\\
&\quad + 4\vp_j \intl^{\vp_j}_0 \left(\sum^N_{k=1} \intl_{\cl F^t 
_{Z_k}(U_k\cap \Omega)} |Z_ku(y)|^2 dV(y)\right) dt \;.
\end{split}
\end{equation}
Note that we adopt the usual convention of denoting by
$C_\vp$ constants which depend only on $\vp$, but whose
actual value may change as the estimates progress.

No point of $\Omega$ is contained in more than $C(C_2)(\vp_j)^{4-4C_{3}}$ of
the sets $\cl F^t_{Z_k}(U_k\cap \Omega)$, $1\le k\le N$, where 
$C(C_2)$ denotes 
a constant depending only on $C_2$. Indeed, if $Q\in\cl F^t_{Z_k}(U_k\cap\Omega)
\cap \cl F^t_{Z_m}(U_m\cap\Omega)$, then by
the 
triangle inequality, the distance between $\cl F^{-t}_{Z_k}(Q)$
and  $\cl F^{-t}_{Z_m}(Q)$ is no more than $2t\le 2\vp_j$ (recall that
$|Z_k| = 1$). Consequently, $B(P_m, \frac{C_2}{10}(\vp_j)^{C_3}) \subseteq 
B(P_k, 
2\vp_j + 2C_2(\vp_j)^{C_3} + \frac{C_2}{10}(\vp_j)^{C_3})$. Since the balls 
$B(P_m, \frac{C_2}{10}(\vp_j)^{C_3})$, $1\le m\le N$, are pairwise disjoint, 
comparison of volumes gives the desired upper bound on how
many of them can be contained in $B(P_k, 2\vp_j +
2C_2(\vp_j)^{C_3} + \frac{C_2}{10}(\vp_j)^{C_3})$.

With this control over the overlap of the sets $\cl F^t_{Z_k}(U_k\cap\Omega)$, 
$1\le k\le N$, for $t$ fixed, we can control the last term
in \eqref{eqinsert}. Denote by $L$ the complex  tangential
field of type (1,0); we may take $L$ to be $C^\infty$ on
$\ovl\Omega$, and so that $|\text{Re } L|$ and $|\text{Im }
L|$ are one in a neighborhood of the boundary that contains
all the sets $U_k\cap \Omega$, $1\le k\le N$. Then, because
$|Z_k|=1$, \begin{equation}\label{eq12} |Z_ku(y)|^2 \le
2|(\text{Re } L) u(y)|^2  + 2|(\text{Im } L)u(y)|^2.
\end{equation} Inserting \eqref{eq12} into \eqref{eqinsert}
and using the estimate on the overlap of the sets $\cl
F^t_{Z_k}(U_k\cap\Omega)$, we obtain

\begin{equation}\label{eq13}
\begin{split}
\intl_{\left(\bigcup\limits^N_{k=1} U_k\right)\cap \Omega}
|u|^2 &\le  2\vp(\|\bar\partial u\|^2_0 +
\|\bar\partial^*u\|^2_0) + C_\vp \|u\|^2_{-1}\\
&\quad + 4\vp_j \intl^{\vp_j}_0 2C(C_2)(\vp_j)^{4-4C_3} \intl_\Omega (|(\text{Re 
} L)u(y)|^2 + |(\text{Im } L)u(y)|^2) dV(y) \ dt\\
& \\
&= 2\vp (\|\bar\partial u\|^2_0 + \|\bar\partial^*u\|^2_0) +
C_\vp \|u\|^2_{-1}\\
& \\
&\quad + 8C(C_2) \vp^{6-4C_3}_j \intl_\Omega (|(\text{Re } L)u(y)|^2 +
|(\text{Im 
} L)u(y)|^2)dV(y).
\end{split}
\end{equation}
Finally, we exploit maximal estimates ([\ref{D1}]):\ since we are in $\bb C^2$, 
the integral in the last term in \eqref{eq13} is dominated by $\|\bar\partial 
u\|^2_0 + \|\bar\partial^*u\|^2_0$. Combining this with \eqref{eq4} and 
\eqref{eq6} shows that there is a constant $C$ independent
of $\vp$ such that for all sufficiently small $\vp>0$, we
have the estimate \begin{equation}\label{eq14}
\|u\|^2_0 \le C(\vp+\vp^{6-4C_3}) (\|\bar\partial u\|^2_0 + \|\bar\partial^*
u\|^2_0) + C_\vp \|u\|^2_{-1}.
\end{equation}
\eqref{eq14} gives the required compactness estimate
\eqref{eq1} (since $C$ does not depend on $\vp$ and
$6-4C_3>0$). This completes the proof of the theorem.
\bigskip
\bigskip
\bigskip

\n {\large References:}

\begin{enumerate}
\item\label{B}
Harold P.\ Boas, Small sets of infinite type are benign for the 
$\bar\partial$-Neumann problem, Proceedings of the American Math.\ Soc. {\bf 
103}, 1988, pp.~569--578.
\item\label{BS}
Harold P.\ Boas and Emil J.\ Straube, Global regularity of the 
$\bar\partial$-Neumann problem:\ a survey of the $\cl L^2$-Sobolev theory, 
{\em Several Complex Variables}, M.\ Schneider and Y.-T.\
Siu, eds., MSRI Publications {\bf 37}, Cambridge University
Press, 1999.
\item\label{C2}
David W. Catlin, Boundary behavior of holomorphic functions on weakly
pseudoconvex domains, Princeton University Ph.D.\ Thesis, 1978.
\item\label{C1}
David W. Catlin, Global regularity of the $\bar\partial$-Neumann problem, 
{\em Complex Analysis of Several Variables}, Y.-T.\ Siu ed.,
Proc.\ Symp.\ Pure Math. {\bf 41}, 1984, pp.~39--49.
\item\label{C3}
David W.\ Catlin, Subelliptic estimates for the $\bar\partial$-Neumann problem 
on pseudoconvex domains, Annals of Mathematics (2) {\bf 126}, 1987, 
pp.~131--191.
\item\label{CS}
So-Chin Chen and Mei-Chi Shaw, {\em Partial Differential Equations in Several 
Complex Variables}, Studies in Advanced Mathematics, American Mathematical 
Society/Interna\-tional Press, 2001.
\item\label{CF}
Michael Christ and Siqi Fu, Compactness in the $\bar\partial$-Neumann problem, 
magnetic Schr\"o\-dinger operators, and the Aharonov--Bohm effect, preprint 
2003.
\item\label{DA}
John P. D'Angelo,  {\em Several Complex Variables and the Geometry of Real 
Hypersurfaces}, Studies in Advanced Mathematics, CRC Press, Boca Raton, FL, 
1993.
\item\label{D1}
M.\ Derridj, Regularit\'e pour $\bar\partial$ dans quelques domaines faiblement 
pseudo-convexes, Journal of Differential Geometry {\bf 13}, 1978, pp.~559--576.
\item\label{FK}
G.B.\ Folland and J.J.\ Kohn, {\em The Neumann Problem for the Cauchy--Riemann 
Complex}, Annals of Mathematics Studies {\bf 75}, Princeton University Press,
1972.
\item\label{FS2}
Siqi Fu and Emil J.\ Straube, Compactness of the $\bar\partial$-Neumann problem
on convex domains, Journal of Functional Analysis {\bf 159}, 1998, pp.~629--641.
\item\label{FS1}
Siqi Fu and Emil J.\ Straube, Compactness in the $\bar\partial$-Neumann problem, 
{\em Complex Analysis and Geometry}, J.\ McNeal ed., Ohio
State Univ.\ Math.\ Research Inst.\ Publ. {\bf 9}, 2001,
pp.~141--160.
\item\label{FS3}
Siqi Fu and Emil J.\ Straube, Semi-classical analysis of Schr\"odinger operators 
and compactness in the $\bar\partial$-Neumann problem, Journal of Math.\ 
Analysis and Application {\bf 271}, 2002, pp.~267--282, correction in Journal of 
Math.\ Analysis and Applications {\bf 280}, 2003, pp.~195--196.
\item\label{HL}
Torsten Hefer and Ingo Lieb, On the compactness of the
$\overline{\partial}$-Neumann operator, Ann. Fac. Sci.
Toulouse Math.(6) {\bf 9}, 2000, pp.~415--432.
\item\label{M} Peter Matheos, A Hartogs domain with no
analytic discs in the boundary for which the
$\bar\partial$-Neumann problem is not compact, University of
California Los Angeles Ph.D.\ Thesis, 1997. \item\label{McN}
Jeffery D.\ McNeal, A sufficient condition for compactness
of the $\bar\partial$-Neumann problem, Journal of Functional
Analysis {\bf 195}, 2002, pp.~190--205.
\item\label{S1}
Nessim Sibony, Une classe de domaines pseudoconvexes, Duke
Math.\ Journal {\bf 55}, 1987, pp.~299--319.
\item\label{St}
Emil J.\ Straube, Plurisubharmonic functions and
subellipticity of the $\overline{\partial}$-Neumann problem 
on non-smooth domains, Mathematical Research Letters {\bf
4}, 1997, pp.~459--467.
\item\label{Z}
William P.\ Ziemer, {\em Weakly Differentiable Functions},
Graduate Texts in Mathematics, vol.\ 120, Springer-Verlag,
1989.

\end{enumerate}

\end{document}